\documentclass[12pt]{amsart}
\usepackage[utf8]{inputenc}
\usepackage[ngerman]{babel}
\usepackage{amsfonts}
\usepackage{latexsym}
\usepackage{amssymb}
\usepackage{amsmath}
\usepackage{float}
\usepackage{tikz}
\usepackage{enumerate}
\usepackage{mathrsfs}
\usepackage{hyperref}
\usepackage[normalem]{ulem}

\usepackage[left=2.5cm, top=2.5cm,bottom=2.5cm,right=2.5cm]{geometry}

\newcommand{\R}{\mathbb R}
\newcommand{\N}{\mathbb N}

\newcommand{\Q}{\mathbb Q}
\newcommand{\Z}{\mathbb Z}

\newtheorem*{satz*}{Satz}
\newtheorem*{folgerung*}{Folgerung}

{
\theoremstyle{definition}

\newtheorem*{rem*}{Bemerkung}

}
\definecolor{sincolor}{rgb}{1,0.83,0.16}
\definecolor{coscolor}{rgb}{0.83,0,0.67}
\definecolor{tancolor}{RGB}{0.0,0,212}
\definecolor{hcolor}{rgb}{0,0.63,0}

\usepackage{nomencl}
\makenomenclature
\begin{document}

%%%%%%%%%%%%%%%%%%%%%%%%%%%%%%%%%%%%%%%%%%%%%5

\title[]{Zur Irrationalität in der Schule}

\author[G. Leobacher]{Gunther Leobacher}
\address{Gunther Leobacher: Institut für Mathematik und Wissenschaftliches Rechnen, Universität Graz, Österreich.}
\email{gunther.leobacher@uni-graz.at}

\author[J. Prochno]{Joscha Prochno}
\address{Joscha Prochno: Institut für Mathematik und Wissenschaftliches Rechnen, Universität Graz, Österreich.} 
\email{joscha.prochno@uni-graz.at}


%\thanks{???}

%\date{\today}

\begin{abstract}
Irrationale Zahlen werden in der Schule bereits in der Sekundarstufe I eingeführt. Allerdings wird typischerweise, mit Ausnahme vielleicht für $\sqrt{2}$, kein mathematischer Beweis zur Irrationalität geführt. Insbesondere wird nicht bewiesen, dass die berühmte Eulersche Zahl $e$ sowie die Kreiszahl $\pi$ irrationale Zahlen sind. In diesem Artikel wollen wir aufzeigen, wie dies mit recht elementaren Methoden der Analysis
möglich ist. Dar\"uber hinaus bieten wir f\"ur viele der analytischen Aussagen geometrische Varianten zur Veranschaulichung, die insbesondere Variabilität im Anspruchsniveau schaffen.
\end{abstract}

\maketitle

\tableofcontents

% % % % % % % % % % % % %
\section{Einleitung}
% % % % % % % % % % % % %

Alles begann bereits im 6.-5. Jahrhundert vor Christus mit den Pythagoreern, die irrationale Zahlen, also solche mit einer nicht periodischen, unendlichen Anzahl von Dezimalstellen, zunächst als \emph{inkommensurable Strecken}, d.h. Strecken mit einem nicht-rationalen Längenverhältnis, entdeckten (etwa beim regelmäßigen Fünfeck). 

Auf den Begriff der irrationalen Zahl treffen Schülerinnen und Schüler heute im Alter von etwa 15 Jahren. Eingeführt wird dieser fast ausschließlich im Zusammenhang mit dem Begriff der Wurzel. Ein typisches Beispiel ist die Irrationalität von $\sqrt{2}$, eine Zahl, die sich auf natürliche Weise aus dem Satz des Pythagoras als Länge der Diagonalen im Quadrat mit Kantenlänge $1$ ergibt. Eine Möglichkeit, sich der Fragestellung nach der Irrationalität von $\sqrt{2}$ zu nähern, ist die der Intervallschachtelung. Aber es bleibt auch eben nur das, eine "`Näherung"'. Selbst eine Vielzahl von Iterationen in der Schachtelung garantiert nämlich nicht, dass es sich nicht doch um einen endlichen oder periodischen Dezimalbruch handelt. Am Ende ist das einzige überzeugende Argument nur der mathematische Beweis. Und für $\sqrt{2}$ ist dieser recht leicht geführt -- und das wird er in Schulen auch oft! Das Vorgehen basiert auf einem Beweis durch Widerspruch, d.h. man nimmt an, dass es Zahlen $a\in\N$ und $b\in\N\setminus\{0\}$ gibt (man beachte hier, dass wir bereits die Positivität ausgenutzt haben), so dass
\[
\sqrt{2} = \frac{a}{b}
\]
gilt. Ohne Beschränkung der Allgemeinheit kann man annehmen, dass nicht sowohl Zähler als auch Nenner gerade Zahlen sind (sonst kürzt man so lange, bis Zähler oder Nenner oder auch beide ungerade sind). Nun nimmt man an, dass $a=2n$ und $b=2m+1$ mit $n,m\in\N$ gilt (die anderen Fälle verlaufen dann analog), d.h.
\[
\sqrt{2} = \frac{2n}{2m+1}.
\]
Quadrieren liefert dann
\[
2 = \frac{4n^2}{4m^2 + 4m +1}
\]
und nach einer elementaren Umformung ergibt sich
\[
4m^2 + 4m +1 = 2n^2.
\]
Auf der linken Seite steht nun aber eine ungerade, auf der rechten eine gerade Zahl. Das ist ein Widerspruch, und $\sqrt{2}$ kann somit nicht rational sein.

Aber wie sieht es etwa mit anderen, gar berühmten Zahlen wie der Eulerschen Zahl $e$ oder der Kreiszahl $\pi$ aus? Auch für viele Schülerinnen und Schüler ist das eine natürliche Frage. Obiger Beweis kann nicht einfach übertragen werden, denn das "`Quadrieren"', oben ein entscheidender Schritt im Beweisargument, der uns eine gerade Zahl liefert, hilft  bei den Zahlen $e$ und $\pi$ eben nicht! Typischerweise wird auf die Irrationalität dieser Zahlen lediglich verwiesen und kein Beweis geführt.
Das ist natürlich auf eine gewisse Art unbefriedigend; der Grund dafür ist wohl, dass man dafür Methoden der höheren Mathematik verwenden muss, die in der Schule nicht zur Verfügung stehen. Oder ist dem nicht so?
Tatsächlich steht ganz am Anfang der Entstehungsgeschichte unseres Artikels der folgende 
Dialog unter den Autoren:
\begin{itemize}
\item[L:] Gestern habe ich während meiner Analysis-Vorlesung für 
LehramtskandidatInnen den Beweis,
dass $\pi$ irrational ist, als Anwendung des Hauptsatzes der Differential-
und Integralrechnung vorgetragen.
\item[P:] Mhm? Nett!
\item[L:] Danach habe ich zu den Studierenden gesagt:
"`Und hier sehen Sie jetzt, was wir schon für tolle Werkzeuge haben.
Das hätten Sie vor wenigen Wochen noch nicht gekonnt!"'
\item[P:] Und, waren Sie begeistert?
\item[L:] Kann ich nicht genau sagen \dots  (Anm.: Die Vorlesung fand Corona bedingt online statt). Aber ich denke mir gerade, das stimmt ja gar nicht, dass sie das vor wenigen Wochen nicht gekonnt hätten!
Der Beweis verwendet nichts, was sie nicht in der Schule auch schon gelernt haben.
\item[P:] Wenn das so ist, dann sollten wir den Beweis für Schüler verständlich aufschreiben.
\end{itemize}
Es soll hier nicht der Eindruck enstehen, dass dieses Aufschreiben
noch nie passiert sei -- 
wir präsentieren mit diesem Artikel lediglich unseren eigenen Versuch.

Wir wollen hier also demonstrieren, dass das Begründen der Irrationalität von $e$ und $\pi$ mit recht elementaren Methoden (zumindest in höheren Klassen) doch möglich und in unseren Augen sinnvoll ist. Des Weiteren zeigt der Beweis für die Irrationalität von $\pi$ auch, wie man sich einer eher zahlentheoretischen Frage mit Methoden der Analysis nähern kann. 

% % % % % % % % % % % % % % % % % % % % % %
\section{Der mathematische Werkzeugkasten}
% % % % % % % % % % % % % % % % % % % % % %

Wir beginnen mit der Wiederholung einiger Grundlagen, die im Beweis zum Einsatz kommen werden. Dazu gehören der binomische Lehrsatz und Eigenschaften des Binomialkoeffizienten, Differentiationseigenschaften der trigonometrischen Funktionen Sinus und Kosinus sowie einige Aussagen zum Integral stetiger Funktionen und der Hauptsatz der Differential- und Integralrechnung. Wir werden stets versuchen auch einen anschaulichen Zugang zu liefern, der, wie wir hoffen, gegebenenfalls auch unabhängig von der Thematik der Irrationalität seine Anwendung in der Schule findet.

Wir verwenden die gängigen Notationen $\N$ für die natürlichen Zahlen (ohne die Null), $\Z$ für die ganzen Zahlen, $\Q$ für die rationalen Zahlen sowie $\R$ für die reellen Zahlen. Des Weiteren verwenden wir, vor allem auch aus Gründen der \"Ubersichtlichkeit, die Notation mit dem Summenzeichen. Für eine Folge $x_1,x_2,\ldots$ reeller Zahlen sowie $m,n\in\N$ mit $m\leq n$ ist also
\[
\sum_{k=m}^n x_k = x_m + x_{m+1} + \dots + x_{n-1}+x_n. 
\]  
Existiert für obige Partialsummen der Grenzwert für $n\to\infty$, so schreibt man kurz
\[
\sum_{k=m}^\infty x_k
\]
für diesen Grenzwert.

% % % % % % % % % % % % % % % % % % %
\subsection{Der binomische Lehrsatz}
% % % % % % % % % % % % % % % % % % %

Wir beginnen unsere Grundlagenreise bei einem mathematischen Satz, der, zumindest im Spezialfall, fast allen Schülerinnen und Schülern ein Begriff ist. Es handelt sich um den binomischen Lehrsatz, den wir für $n=2$ auch als (erste) "`binomische Formel"' kennen. Der Satz sagt uns, wie wir die $n$-te Potenz einer Summe zweier Zahlen als ein Polynom $n$-ten Grades schreiben können. Eine zentrale Rolle in dieser Darstellung spielen die Binomialkoeffizienten. Auf einen Beweis, der mittels vollständiger Induktion geführt werden kann, verzichten wir hier. 
 
\begin{satz*}
Seien $a,b \in\R$ und $n\in \N\cup\{0\}$. Dann gilt  
\[
(a+b)^n=a^n+na^{n-1}b+\ldots +{n \choose k} a^{n-k}b^k+\ldots+nab^{n-1}+b^n = \sum_{k=0}^n{n\choose k}a^{n-k}b^k\;,
\]
wobei die Binomialkoeffizienten 
${n \choose k}:=\frac{n!}{(n-k)!k!}$  die Rekursion 
\begin{equation}\label{eq:binom}
{n \choose k}={n-1 \choose k-1}+{n-1 \choose k}
\end{equation}
mit den Startwerten ${n \choose 0}=1$ und ${n\choose n}=1$ erfüllt.
\end{satz*}

Im Fall, dass $n=2$ ist und die Zahlen $a,b$ positiv sind, kann man den binomischen Lehrsatz leicht geometrisch visualisieren. Dazu nehmen wir uns ein Quadrat der Seitenlänge $a+b$. Dann sehen wir, dass wir dieses zerlegen können in ein Quadrat der Seitenlänge $a$, ein Quadrat der Seitenlänge $b$ sowie \underline{zwei} Rechtecke mit Seiten $a$ und $b$.
\smallskip
\begin{center}
%% Creator: Inkscape 1.0 (4035a4f, 2020-05-01), www.inkscape.org
%% PDF/EPS/PS + LaTeX output extension by Johan Engelen, 2010
%% Accompanies image file 'binom2.pdf' (pdf, eps, ps)
%%
%% To include the image in your LaTeX document, write
%%   \input{<filename>.pdf_tex}
%%  instead of
%%   \includegraphics{<filename>.pdf}
%% To scale the image, write
%%   \def\svgwidth{<desired width>}
%%   \input{<filename>.pdf_tex}
%%  instead of
%%   \includegraphics[width=<desired width>]{<filename>.pdf}
%%
%% Images with a different path to the parent latex file can
%% be accessed with the `import' package (which may need to be
%% installed) using
%%   \usepackage{import}
%% in the preamble, and then including the image with
%%   \import{<path to file>}{<filename>.pdf_tex}
%% Alternatively, one can specify
%%   \graphicspath{{<path to file>/}}
%% 
%% For more information, please see info/svg-inkscape on CTAN:
%%   http://tug.ctan.org/tex-archive/info/svg-inkscape
%%
\begingroup%
  \makeatletter%
  \providecommand\color[2][]{%
    \errmessage{(Inkscape) Color is used for the text in Inkscape, but the package 'color.sty' is not loaded}%
    \renewcommand\color[2][]{}%
  }%
  \providecommand\transparent[1]{%
    \errmessage{(Inkscape) Transparency is used (non-zero) for the text in Inkscape, but the package 'transparent.sty' is not loaded}%
    \renewcommand\transparent[1]{}%
  }%
  \providecommand\rotatebox[2]{#2}%
  \newcommand*\fsize{\dimexpr\f@size pt\relax}%
  \newcommand*\lineheight[1]{\fontsize{\fsize}{#1\fsize}\selectfont}%
  \ifx\svgwidth\undefined%
    \setlength{\unitlength}{174.24301065bp}%
    \ifx\svgscale\undefined%
      \relax%
    \else%
      \setlength{\unitlength}{\unitlength * \real{\svgscale}}%
    \fi%
  \else%
    \setlength{\unitlength}{\svgwidth}%
  \fi%
  \global\let\svgwidth\undefined%
  \global\let\svgscale\undefined%
  \makeatother%
  \begin{picture}(1,0.90333039)%
    \lineheight{1}%
    \setlength\tabcolsep{0pt}%
    \put(0,0){\includegraphics[width=\unitlength,page=1]{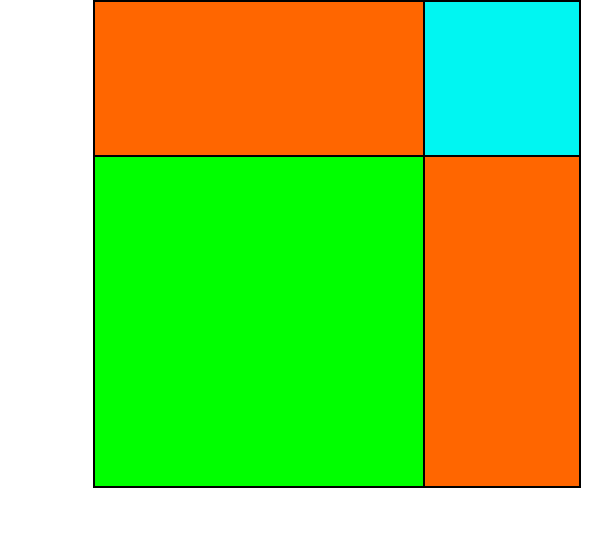}}%
    \put(0.11802221,0.36313313){\makebox(0,0)[rt]{\lineheight{1.25}\smash{\begin{tabular}[t]{r}$a$\end{tabular}}}}%
    \put(0.11494769,0.75913168){\makebox(0,0)[rt]{\lineheight{1.25}\smash{\begin{tabular}[t]{r}$b$\end{tabular}}}}%
    \put(0.43667735,0.01481796){\makebox(0,0)[t]{\lineheight{1.25}\smash{\begin{tabular}[t]{c}$a$\end{tabular}}}}%
    \put(0.84066957,0.01112795){\makebox(0,0)[t]{\lineheight{1.25}\smash{\begin{tabular}[t]{c}$b$\end{tabular}}}}%
    \put(0.44744158,0.35977947){\makebox(0,0)[t]{\lineheight{1.25}\smash{\begin{tabular}[t]{c}$a^2$\end{tabular}}}}%
    \put(0.83544633,0.76192706){\makebox(0,0)[t]{\lineheight{1.25}\smash{\begin{tabular}[t]{c}$b^2$\end{tabular}}}}%
    \put(0.43329874,0.76192699){\makebox(0,0)[t]{\lineheight{1.25}\smash{\begin{tabular}[t]{c}$a\cdot b$\end{tabular}}}}%
    \put(0.83237181,0.36592838){\makebox(0,0)[t]{\lineheight{1.25}\smash{\begin{tabular}[t]{c}$a\cdot b$\end{tabular}}}}%
  \end{picture}%
\endgroup%

\end{center}

Die Gleichung \eqref{eq:binom} im binomischen Lehrsatz wird gewöhnlich durch das Pascal'sche Dreieck
illustriert:
\begin{center}
\begin{tabular}{c|ccccccccc}
$n\backslash k$&0&1&2&3&4&5&$\cdots$ \\\hline
0 & 1 &   &    &   &   &   &  \\
1 & 1 & 1 &    &   &   &   &  \\
2 & 1 & 2 &  1 &   &   &   &  \\
3 & 1 & 3 &  3 & 1 &   &   &  \\
4 & 1 & 4 &  6 & 4 & 1 &   &  \\
5 & 1 & 5 & 10 & 10 & 5 & 1 & \\
$\vdots$ &$\vdots$ &$\vdots$ &$\vdots$ &$\vdots$ &$\vdots$ &$\vdots$ &$\ddots$\\
\end{tabular}
\end{center}

\bigskip

Der in den Binomialkoeffizienten auftauchende Ausdruck $n!$ bezeichnet die {\em Faktorielle von $n$}, auch 
{\em $n$ Fakulät}, und ist definiert als das Produkt der ersten 
$n$ natürlichen Zahlen, also $n!=1\cdot 2\cdot \ldots\cdot (n-1)\cdot n\;$. Für die Zahl 0 definiert man zusätzlich $0!=1\,$. Die ersten paar Faktoriellen werden in der nachstehenden Tabelle anderen Folgen gegenüber gestellt.
\begin{center}
\begin{tabular}{r|r r r r r }
$n$&$n^2$ & $n^3$ & $2^n$ &$n!$\\\hline
0 &  0 &    0 &   1 &    1 \\
1 &  1 &    1 &   2 &    1 \\
2 &  4 &    8 &   4 &    2 \\
3 &  9 &   27 &   8 &    6 \\
4 & 16 &   64 &  16 &   24 \\
5 & 25 &  125 &  32 &  120 \\
6 & 36 &  216 &  64 &  720 \\
7 & 49 &  343 & 128 & 5040 \\
\end{tabular}
\end{center}

\bigskip

Man sieht, dass die Faktorielle sehr schnell wächst, und zwar, wie der
folgende Satz zeigt, schneller als eine Exponentialfolge mit beliebiger Basis.

\begin{satz*}
Für beliebige Zahlen $a,c\in(0;\infty)$ gibt es stets eine natürliche
Zahl $n\in\N$ mit \[c\cdot a^n< n!\;. \]
\end{satz*}

\begin{proof}
Es genügt zu zeigen, dass $a^n/n!$ für $n\to\infty$ eine Nullfolge bildet. Denn dann existiert ein $n_0\in\N$, so dass für jedes $n\geq n_0$, $a^n/n!\leq 1/(2c)$ gilt. Wir wählen zunächst eine Zahl $n_0\in\N$, die größer als $2 a$ ist, sodass also $a/n_0 < 1/2$. Für alle $n\in\N$ mit $n\geq n_0$ gilt dann
\[
\frac{a^n}{n!} = \frac{a^{n_0}}{n_0!}\cdot\frac{a^{n-n_0}}{(n_0+1)\cdot\ldots\cdot n} < \frac{a^{n_0}}{n_0!} \cdot\Big(\frac{1}{2}\Big)^{n-n_0}.
\]
Für $n\to\infty$ ist letzerer Ausdruck eine feste Zahl multipliziert mit einer Nullfolge, d.h. insgesamt eine Nullfolge.
\end{proof}

Eine wichtige Erkenntnis aus \eqref{eq:binom} ist aber auch die, dass
die Binomialkoeffizienten stets ganze Zahlen sind. Dass für Zahlen $k,n\in\N$ mit $0\le k\le n$ der Bruch 
$\frac{n!}{k!}$ stets eine ganze Zahl ist, folgt unmittelbar aus der
Definition der Faktoriellen. Die Ganzzahligkeit der Binomialkoeffizienten,
für welche eben durch $(n-k)!$ auch noch dividiert wird, erfordert eine 
zusätzliche Überlegung, wie eben die obige.

\begin{folgerung*}
Die Binomialkoeffizienten ${n\choose k}$ sind für alle $k\in\N$ mit $0\le k\le n$ 
ganzzahlig.
\end{folgerung*}

Eine weitere wichtige Konsequenz aus dem binomischen Lehrsatz ist die Ableitungsformel für Potenzfunktionen.

\begin{folgerung*}
Für jedes $n\in \N$ ist die Ableitung der Potenzfunktion gegeben durch 
\[
(x^n)'=nx^{n-1}\;.
\]
Allgemein gilt für die $\ell$-te Ableitung der Potenzfunktion 
\begin{equation}\label{eq:ableitung-potenz}
(x^n)^{(\ell)}=\begin{cases}
\frac{n!}{(n-\ell)!}x^{n-\ell}&:\, 0\le \ell\le n\\
0 &:\, \ell> n \;.
\end{cases}
\end{equation}
\end{folgerung*}

\begin{proof} Für $n=0$ und $n=1$ rechnet man $1'=0$  bzw $x'=1$  ganz einfach nach.
Für $n\ge 2$ fixieren wir ein $x\in\R$ und betrachten den Differenzenquotienten $\frac{(x+h)^n-x^n}{h}$.
Für dessen Zähler gilt mithilfe des binomischen Lehrsatzes
\[
(x+h)^n-x^n
={n\choose 1} x^{n-1}h+{n\choose 2} x^{n-2}h^2+\ldots+ {n\choose n} h^n\;,
\]
so dass 
\[
\frac{(x+h)^n-x^n}{h}={n\choose 1} x^{n-1}+h\bigg({n\choose 2} x^{n-2}+\ldots+ {n\choose n} h^{n-2}\bigg)\;.
\]
Somit gilt 
\[
(x^n)'=\lim_{h\to 0}\frac{(x+h)^n-x^n}{h}={n\choose 1} x^{n-1}=nx^{n-1}\;,
\]
wobei wir ${n\choose 1}=n$ aus dem Pascal'schen Dreieck ablesen.

Wenn wir nun noch einmal ableiten, so erhalten wir 
\[
(x^n)''=\begin{cases}
n (n-1)x^{n-2}&:\, n\ge 2\\
0 &:\,  n<2 \;,
\end{cases}
\] 
wobei wir beobachten, dass $n (n-1)=\frac{n!}{(n-2)!}\,$.
Ableiten von $x^{(\ell-1)}=\frac{n!}{(n-\ell-1)!}x^{n-\ell+1}$ liefert nun
\[
x^{(\ell)}=\big(x^{(\ell-1)})'
=\begin{cases}
\frac{n!}{(n-\ell+1)!}(n-\ell+1)x^{n-\ell}&:\, n-\ell\ge 0\\
0&:\, n-\ell< 0\;.
\end{cases}
\]
Nun ist aber $\frac{n!}{(n-\ell+1)!}(n-\ell+1)=\frac{n!}{(n-\ell)!}$, und die Behauptung ist gezeigt.
\end{proof}

% % % % % % % % % % % % % % % % % % % % % % % % % %
\subsection{Sinus, Cosinus und die Kreiszahl $\pi$}
% % % % % % % % % % % % % % % % % % % % % % % % % %

Die unumstritten bekanntesten trigonometrischen Funktionen sind die Sinusfunktion und die Kosinusfunktion. In der Schule werden diese Funktionen geometrisch im Zusammenhang mit rechtwinkligen Dreiecken eingeführt. So ist etwa der Sinus eines Winkels gegeben durch das Verhältnis der Länge der Gegenkathete (also die Kathete, die dem Winkel gegenüberliegt) zur Länge der Hypotenuse (die Seite gegenüber dem rechten Winkel). Eine beliebte Veranschaulichung von Sinus und Kosinus ist die am Einheitskreis.

\begin{center}
%% Creator: Inkscape 1.0 (4035a4f, 2020-05-01), www.inkscape.org
%% PDF/EPS/PS + LaTeX output extension by Johan Engelen, 2010
%% Accompanies image file 'sinuscosinus.pdf' (pdf, eps, ps)
%%
%% To include the image in your LaTeX document, write
%%   \input{<filename>.pdf_tex}
%%  instead of
%%   \includegraphics{<filename>.pdf}
%% To scale the image, write
%%   \def\svgwidth{<desired width>}
%%   \input{<filename>.pdf_tex}
%%  instead of
%%   \includegraphics[width=<desired width>]{<filename>.pdf}
%%
%% Images with a different path to the parent latex file can
%% be accessed with the `import' package (which may need to be
%% installed) using
%%   \usepackage{import}
%% in the preamble, and then including the image with
%%   \import{<path to file>}{<filename>.pdf_tex}
%% Alternatively, one can specify
%%   \graphicspath{{<path to file>/}}
%% 
%% For more information, please see info/svg-inkscape on CTAN:
%%   http://tug.ctan.org/tex-archive/info/svg-inkscape
%%
\begingroup%
  \makeatletter%
  \providecommand\color[2][]{%
    \errmessage{(Inkscape) Color is used for the text in Inkscape, but the package 'color.sty' is not loaded}%
    \renewcommand\color[2][]{}%
  }%
  \providecommand\transparent[1]{%
    \errmessage{(Inkscape) Transparency is used (non-zero) for the text in Inkscape, but the package 'transparent.sty' is not loaded}%
    \renewcommand\transparent[1]{}%
  }%
  \providecommand\rotatebox[2]{#2}%
  \newcommand*\fsize{\dimexpr\f@size pt\relax}%
  \newcommand*\lineheight[1]{\fontsize{\fsize}{#1\fsize}\selectfont}%
  \ifx\svgwidth\undefined%
    \setlength{\unitlength}{220.0482815bp}%
    \ifx\svgscale\undefined%
      \relax%
    \else%
      \setlength{\unitlength}{\unitlength * \real{\svgscale}}%
    \fi%
  \else%
    \setlength{\unitlength}{\svgwidth}%
  \fi%
  \global\let\svgwidth\undefined%
  \global\let\svgscale\undefined%
  \makeatother%
  \begin{picture}(1,0.94801156)%
    \lineheight{1}%
    \setlength\tabcolsep{0pt}%
    \put(0,0){\includegraphics[width=\unitlength,page=1]{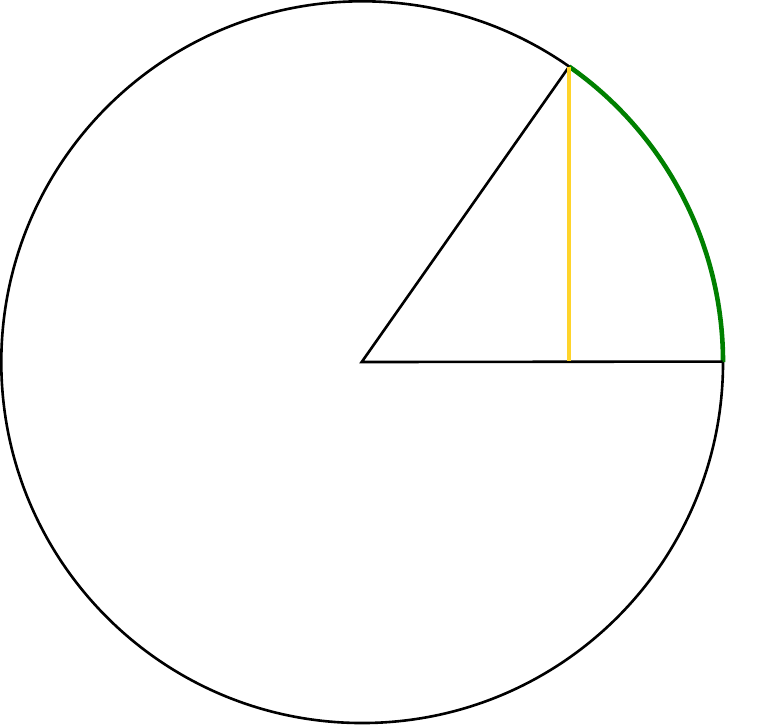}}%
    \put(0.75957614,0.59208043){\color[rgb]{1,0.86666667,0.33333333}\makebox(0,0)[lt]{\lineheight{1.25}\smash{\begin{tabular}[t]{l}$\sin(x)$\end{tabular}}}}%
    \put(0.60181864,0.41801158){\color[rgb]{0.83137255,0,0.66666667}\makebox(0,0)[t]{\lineheight{1.25}\smash{\begin{tabular}[t]{c}$\cos(x)$\end{tabular}}}}%
    \put(0.9102735,0.71332014){\color[rgb]{0,0.50196078,0}\makebox(0,0)[lt]{\lineheight{1.25}\smash{\begin{tabular}[t]{l}$x$\end{tabular}}}}%
    \put(0,0){\includegraphics[width=\unitlength,page=2]{sinuscosinus.pdf}}%
  \end{picture}%
\endgroup%

\end{center}

In der (höheren) Analysis führt man diese Funktionen über eine Reihendarstellung ein, genauer über die sogenannten Potenzreihen. Die Darstellung ist dann
\[
\sin(x) = \sum_{k=0}^\infty (-1)^k\frac{x^{2k+1}}{(2k+1)!} = x - \frac{x^3}{3!} + \frac{x^5}{5!}\mp\dots
\]
für den Sinus und
\[
\cos(x) = \sum_{k=0}^\infty (-1)^k\frac{x^{2k}}{(2k)!} = 1 - \frac{x^2}{2!} + \frac{x^4}{4!}\mp\dots
\]
für den Kosinus einer Zahl $x\in\R$. Was wir im Rahmen dieses Artikels benötigen, sind die Ableitungen von Sinus und Kosinus, die auch in den entsprechenden Klassenstufen an Schulen behandelt werden. Es gelten
\[
\sin'(x) = \cos(x) \qquad\text{und}\qquad \cos'(x)= - \sin(x).
\]
Selbst der entsprechende Beweis kann mit Hilfe der Definition am Einheitskreis und einem Additionstheorem geführt werden. Wir skizzieren dies hier lediglich im Falle der Sinusfunktion. Dazu erinnern wir zunächst an ein Additionstheorem für den Sinus. Es gilt für alle $a,b\in\R$, dass  
\[
\sin(a+b) = \sin(a)\cos(b) + \sin(b)\cos(a).
\]
Damit erhält man nun für alle $x\in\R$ und $h\in\R\setminus\{0\}$, dass
\[
\frac{\sin(x+h)-\sin(x)}{h} = \frac{\sin(x)\cos(h) + \sin(h)\cos(x) - \sin(x)}{h} = \sin(x)\frac{\cos(h)-1}{h} + \cos(x)\frac{\sin(h)}{h}.
\]
Betrachtet man nun den Grenzwert für $h$ gegen $0$ und nutzt die Tatsache, dass
\begin{equation}\label{eq:limits sine and consine}
\lim_{h\to 0}\frac{\cos(h)-1}{h} = 0 \qquad\text{und}\qquad \lim_{h\to 0}\frac{\sin(h)}{h} = 1
\end{equation}
gilt (das kann man sich natürlich auch überlegen -- etwa mit Hilfe von Abschätzungen, die man an den Reihendarstellungen "`ablesen"' kann), 
so erhält man
\[
\sin'(x) = \lim_{h\to 0} \frac{\sin(x+h)-\sin(x)}{h} = \cos(x).
\]

Sicherlich ist obiges Argument nur in höheren Klassen durchführbar oder vielleicht auch nur an mathematisch-technisch ausgerichteten Schulen. Man kann sich aber diese Grenzwerte und Ableitungen auch geometrisch am Einheitskreis plausibel machen und somit diesen Sachverhalt für Schülerinnen und Schüler veranschaulichen. 

\begin{center}
%% Creator: Inkscape 1.0 (4035a4f, 2020-05-01), www.inkscape.org
%% PDF/EPS/PS + LaTeX output extension by Johan Engelen, 2010
%% Accompanies image file 'sinusxdurchx.pdf' (pdf, eps, ps)
%%
%% To include the image in your LaTeX document, write
%%   \input{<filename>.pdf_tex}
%%  instead of
%%   \includegraphics{<filename>.pdf}
%% To scale the image, write
%%   \def\svgwidth{<desired width>}
%%   \input{<filename>.pdf_tex}
%%  instead of
%%   \includegraphics[width=<desired width>]{<filename>.pdf}
%%
%% Images with a different path to the parent latex file can
%% be accessed with the `import' package (which may need to be
%% installed) using
%%   \usepackage{import}
%% in the preamble, and then including the image with
%%   \import{<path to file>}{<filename>.pdf_tex}
%% Alternatively, one can specify
%%   \graphicspath{{<path to file>/}}
%% 
%% For more information, please see info/svg-inkscape on CTAN:
%%   http://tug.ctan.org/tex-archive/info/svg-inkscape
%%
\begingroup%
  \makeatletter%
  \providecommand\color[2][]{%
    \errmessage{(Inkscape) Color is used for the text in Inkscape, but the package 'color.sty' is not loaded}%
    \renewcommand\color[2][]{}%
  }%
  \providecommand\transparent[1]{%
    \errmessage{(Inkscape) Transparency is used (non-zero) for the text in Inkscape, but the package 'transparent.sty' is not loaded}%
    \renewcommand\transparent[1]{}%
  }%
  \providecommand\rotatebox[2]{#2}%
  \newcommand*\fsize{\dimexpr\f@size pt\relax}%
  \newcommand*\lineheight[1]{\fontsize{\fsize}{#1\fsize}\selectfont}%
  \ifx\svgwidth\undefined%
    \setlength{\unitlength}{192.3166375bp}%
    \ifx\svgscale\undefined%
      \relax%
    \else%
      \setlength{\unitlength}{\unitlength * \real{\svgscale}}%
    \fi%
  \else%
    \setlength{\unitlength}{\svgwidth}%
  \fi%
  \global\let\svgwidth\undefined%
  \global\let\svgscale\undefined%
  \makeatother%
  \begin{picture}(1,1.070212)%
    \lineheight{1}%
    \setlength\tabcolsep{0pt}%
    \put(0,0){\includegraphics[width=\unitlength,page=1]{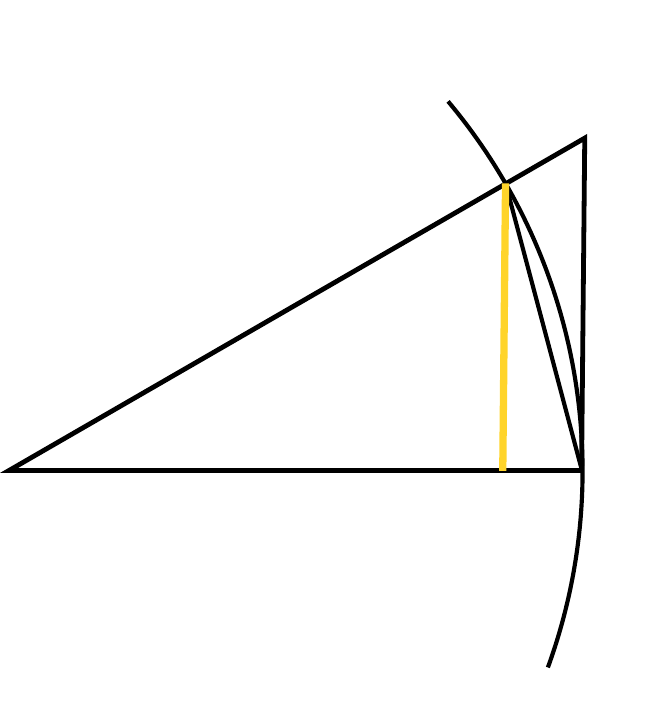}}%
    \put(0.7359454,0.53867081){\color[rgb]{1,0.86666667,0.33333333}\makebox(0,0)[rt]{\lineheight{1.25}\smash{\begin{tabular}[t]{r}$\sin(h)$\end{tabular}}}}%
    \put(0.40605008,0.27863538){\color[rgb]{0.83137255,0,0.66666667}\makebox(0,0)[t]{\lineheight{1.25}\smash{\begin{tabular}[t]{c}$\cos(h)$\end{tabular}}}}%
    \put(0.89368559,0.46581123){\color[rgb]{0,0.50196078,0}\makebox(0,0)[lt]{\lineheight{1.25}\smash{\begin{tabular}[t]{l}$h$\end{tabular}}}}%
    \put(0,0){\includegraphics[width=\unitlength,page=2]{sinusxdurchx.pdf}}%
    \put(0.75102853,0.28643458){\makebox(0,0)[t]{\lineheight{1.25}\smash{\begin{tabular}[t]{c}$C$\end{tabular}}}}%
    \put(0.89169652,0.34028225){\makebox(0,0)[lt]{\lineheight{1.25}\smash{\begin{tabular}[t]{l}$A$\end{tabular}}}}%
    \put(0.73393316,0.78627981){\makebox(0,0)[rt]{\lineheight{1.25}\smash{\begin{tabular}[t]{r}$B$\end{tabular}}}}%
    \put(0.91252083,0.88471385){\makebox(0,0)[t]{\lineheight{1.25}\smash{\begin{tabular}[t]{c}$D$\end{tabular}}}}%
    \put(0.89498915,0.58546654){\makebox(0,0)[lt]{\lineheight{1.25}\smash{\begin{tabular}[t]{l}$E$\end{tabular}}}}%
    \put(0,0){\includegraphics[width=\unitlength,page=3]{sinusxdurchx.pdf}}%
  \end{picture}%
\endgroup%

\end{center}

Wir wollen uns nun geometrisch überlegen, dass die zweite Gleichung in \eqref{eq:limits sine and consine} gilt. In unserer Argumentation folgen wir Stewart \cite{St2015}. Wir wollen die Länge der kürzesten Strecke zwischen zwei Punkten $P_1$ und $P_2$ in der Ebene mit $|P_1P_2|$ bezeichnen.  Zunächst nehmen wir an, dass $h\in(0,\pi/2)$ gilt ($h$ messen wir im Bogenmaß). Wir können der Grafik unmittelbar entnehmen, dass $|BC| = \sin(h)\cdot 1 = \sin(h)$ ist. Nun ist aber $|BC|$ kleiner als $|BA|$, was kleiner als der Bogen von $A$ nach $B$, d.h. $h$ ist. Das liefert sofort die Ungleichung $\sin(h) < h$ und somit
\begin{equation}\label{ineq: sin h durch h}
\frac{\color{sincolor}\sin(h)}{\color{hcolor}h}<1.
\end{equation} 
Den Punkt, an dem sich die Tagenten am Kreis durch $A$ und $B$ treffen, wollen wir mit $E$ bezeichnen. Dann gilt also, dass $h$ kleiner ist als $|AE| + |BE|$. Es folgt weiter
\[
{\color{hcolor}h} < |AE| + |BE| < |AE| + |ED| = {\color{tancolor}\tan(h)}\,. 
\]
Also gilt, wegen der Identität $\tan(x) = \sin(x)/\cos(x)$ für $x\in\R$, dass 
\[
{\color{hcolor}h} < \frac{\color{sincolor}\sin(h)}{\color{coscolor}\cos(h)}.
\]
Kombinieren wir dies mit \eqref{ineq: sin h durch h}, so ergibt sich
\[
{\color{coscolor}\cos(h)} < \frac{\color{sincolor}\sin(h)}{\color{hcolor}h} < 1.
\]
Nun ergibt sich wegen $\lim_{h\to0}\cos(h)=1$ sowie $\lim_{h\to0}1 = 1$, dass
\[
\lim_{h\downarrow 0} \frac{\color{sincolor}\sin(h)}{\color{hcolor}h} = 1
\]
gelten muss, denn $\sin(h)/h$ wird für $h\downarrow 0$ \footnote{Hier bezeichnet $\lim_{h\downarrow 0}$ den Grenzwert von rechts, d.h. nur für positive $h$.} von zwei Größen eingeschachtelt, die beiden gegen $1$ konvergieren. Weil die Funktion $\sin(h)/h$ aber eine gerade Funktion ist, stimmt der linksseitige Grenzwert mit dem rechtsseitigen überein, es gilt also
\[
\lim_{h\to 0} \frac{\sin(h)}{h} = 1.
\]
Die erste Gleichung in \eqref{eq:limits sine and consine} folgt jetzt 
aus dem Additionstheorem $\sin^2(x)+\cos^2(x)=1$, $x\in\R$.
Man kann sie aber auch wie folgt direkt sehen: Es ist wie oben
$h\ge |AB|$, also auch $h^2\ge |AB|^2$, und andererseits
mit dem Satz von Pythagoras (im ersten und im vorletzten "`$=$"') 
\begin{align*}
|AB|^2=|BC|^2+|AC|^2&=\sin^2(h)+\big(1-\cos(h)\big)^2\\
&= \sin^2(h)+1-2 \cos(h)+\cos^2(h)\\
&=2\cdot\big(1-\cos(h)\big)\,,
\end{align*}
so dass $h^2\ge 2 \cdot\big(1-\cos(h)\big)\ge 0$ und damit 
\[
\frac{h}{2}\ge \frac{1-\cos(h)}{h}\ge 0\;.
\]  

Die nächste Grafik unten zeigt, wie man sich die Ableitung des Sinus geometrisch  überlegen kann, wobei dies am Einheitskreis unter Benutzung eines Strahlensatzes geschieht. 

\begin{center}
%% Creator: Inkscape 1.0 (4035a4f, 2020-05-01), www.inkscape.org
%% PDF/EPS/PS + LaTeX output extension by Johan Engelen, 2010
%% Accompanies image file 'ableitung-sinus-gr.pdf' (pdf, eps, ps)
%%
%% To include the image in your LaTeX document, write
%%   \input{<filename>.pdf_tex}
%%  instead of
%%   \includegraphics{<filename>.pdf}
%% To scale the image, write
%%   \def\svgwidth{<desired width>}
%%   \input{<filename>.pdf_tex}
%%  instead of
%%   \includegraphics[width=<desired width>]{<filename>.pdf}
%%
%% Images with a different path to the parent latex file can
%% be accessed with the `import' package (which may need to be
%% installed) using
%%   \usepackage{import}
%% in the preamble, and then including the image with
%%   \import{<path to file>}{<filename>.pdf_tex}
%% Alternatively, one can specify
%%   \graphicspath{{<path to file>/}}
%% 
%% For more information, please see info/svg-inkscape on CTAN:
%%   http://tug.ctan.org/tex-archive/info/svg-inkscape
%%
\begingroup%
  \makeatletter%
  \providecommand\color[2][]{%
    \errmessage{(Inkscape) Color is used for the text in Inkscape, but the package 'color.sty' is not loaded}%
    \renewcommand\color[2][]{}%
  }%
  \providecommand\transparent[1]{%
    \errmessage{(Inkscape) Transparency is used (non-zero) for the text in Inkscape, but the package 'transparent.sty' is not loaded}%
    \renewcommand\transparent[1]{}%
  }%
  \providecommand\rotatebox[2]{#2}%
  \newcommand*\fsize{\dimexpr\f@size pt\relax}%
  \newcommand*\lineheight[1]{\fontsize{\fsize}{#1\fsize}\selectfont}%
  \ifx\svgwidth\undefined%
    \setlength{\unitlength}{302.00560164bp}%
    \ifx\svgscale\undefined%
      \relax%
    \else%
      \setlength{\unitlength}{\unitlength * \real{\svgscale}}%
    \fi%
  \else%
    \setlength{\unitlength}{\svgwidth}%
  \fi%
  \global\let\svgwidth\undefined%
  \global\let\svgscale\undefined%
  \makeatother%
  \begin{picture}(1,0.62007302)%
    \lineheight{1}%
    \setlength\tabcolsep{0pt}%
    \put(0,0){\includegraphics[width=\unitlength,page=1]{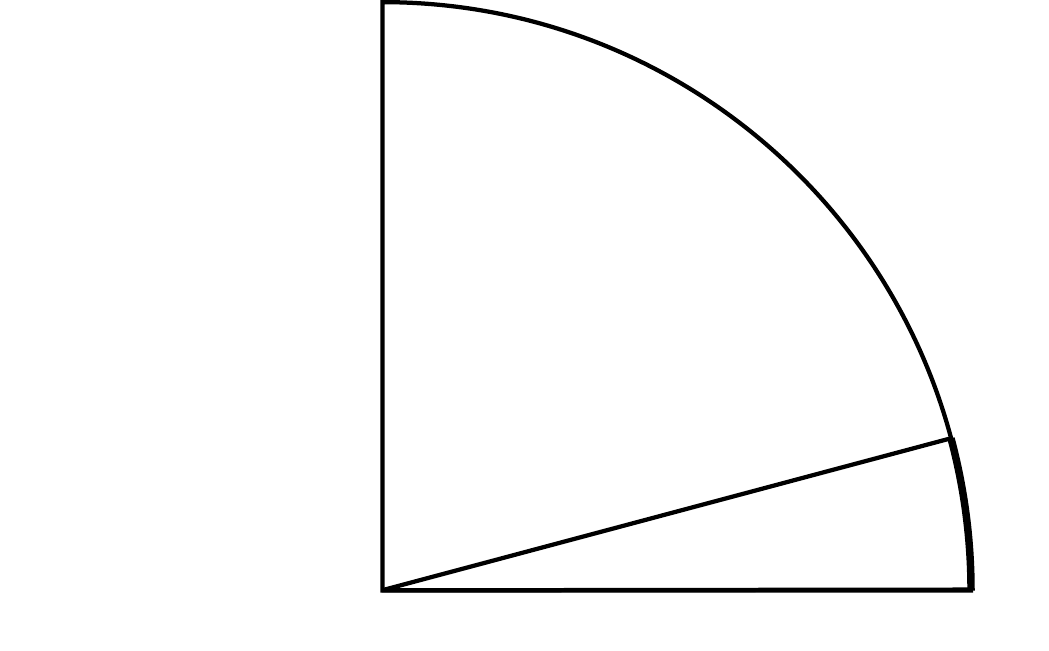}}%
    \put(0.6370177,0.00732732){\color[rgb]{0.83137255,0,0.66666667}\makebox(0,0)[t]{\lineheight{1.25}\smash{\begin{tabular}[t]{c}$\approx\cos(x+\frac{h}{2})$\end{tabular}}}}%
    \put(0.93462326,0.11372866){\makebox(0,0)[lt]{\lineheight{1.25}\smash{\begin{tabular}[t]{l}$x$\end{tabular}}}}%
    \put(0,0){\includegraphics[width=\unitlength,page=2]{ableitung-sinus-gr.pdf}}%
    \put(0.92366406,0.19674688){\color[rgb]{0,0,0}\makebox(0,0)[lt]{\lineheight{1.25}\smash{\begin{tabular}[t]{l}$A$\end{tabular}}}}%
    \put(0.79901714,0.42364448){\color[rgb]{0,0,0}\makebox(0,0)[lt]{\lineheight{1.25}\smash{\begin{tabular}[t]{l}$B$\end{tabular}}}}%
    \put(0,0){\includegraphics[width=\unitlength,page=3]{ableitung-sinus-gr.pdf}}%
    \put(0.35976809,0.30494957){\color[rgb]{1,0.8,0}\makebox(0,0)[rt]{\lineheight{1.25}\smash{\begin{tabular}[t]{r}$\sin(x+h)-\sin(x)$\end{tabular}}}}%
    \put(0.8803382,0.31346851){\color[rgb]{0,0.62745098,0}\makebox(0,0)[lt]{\lineheight{1.25}\smash{\begin{tabular}[t]{l}$h$\end{tabular}}}}%
    \put(0.78398125,0.21134158){\color[rgb]{0,0,0}\makebox(0,0)[rt]{\lineheight{1.25}\smash{\begin{tabular}[t]{r}$C$\end{tabular}}}}%
    \put(0.35335108,0.03593525){\color[rgb]{0,0,0}\makebox(0,0)[rt]{\lineheight{1.25}\smash{\begin{tabular}[t]{r}$O$\end{tabular}}}}%
    \put(0.85788588,0.06806065){\color[rgb]{0,0,0}\makebox(0,0)[lt]{\lineheight{1.25}\smash{\begin{tabular}[t]{l}$P$\end{tabular}}}}%
    \put(0.84335348,0.2632562){\color[rgb]{0,0,0}\makebox(0,0)[rt]{\lineheight{1.25}\smash{\begin{tabular}[t]{r}$Q$\end{tabular}}}}%
    \put(0,0){\includegraphics[width=\unitlength,page=4]{ableitung-sinus-gr.pdf}}%
  \end{picture}%
\endgroup%

\end{center}

Wir nutzen hier die Ähnlichkeit der zwei rechtwinkeligen
Dreiecke $\Delta(ABC)$ und $\Delta(OPQ)$. Diese liefert uns mittels des Strahlensatzes, dass sich $\sin(x+h)-\sin(x)$ zu $|AB|$ verhält, wie sich $|OP|$ zu $|OQ| = 1$ verhält. Wir erhalten nun
\begin{align*}
\frac{\color{sincolor}\sin(x+h)-\sin(x)}{\color{hcolor}h}
\approx 
\frac{\color{sincolor}\sin(x+h)-\sin(x)}{|AB|}
=\frac{|OP|}{|OQ|}
\approx\frac{\color{coscolor}\cos\big(x+\frac{h}{2}\big)}{1}.
\end{align*}

Wir kommen nun zu der vielleicht berühmtesten mathematischen Konstanten, der Kreiszahl $\pi$. Eine Möglichkeit der Definition dieser Zahl ist geometrisch, nämlich als das Verhältnis des Umfangs eines Kreises zu seinem Durchmesser, d.h.
\[
\pi := \frac{\text{Umfang des Kreises}}{\text{Durchmesser des Kreises}}. 
\]
Ein möglicher analytischer Weg zur Definition führt wieder über die trigonometrischen Funktionen und Nullstellen. Hat man etwa die Kosinusfunktion über die Potenzreihe dargestellt, dann kann man $\pi$ als das Doppelte der kleinsten positiven Nullstelle des Kosinus festlegen. Aber welchen Wert hat diese Zahl? Erste Näherungen finden sich bereits im berühmten Papyrus Rhind aus dem 16. Jahrhundert vor Christus, nämlich 
\[
\Big(\frac{16}{9}\Big)^2 \approx 3{,}1605.
\]
Ist $\pi$ also vielleicht wirklich eine rationale Zahl, aber eben nicht die oben gegebene? Archimedes begegnete der Kreiszahl $\pi$ durch einer Approximation des Kreises durch Vielecke \cite{AH2001} und erhielt die Näherung
\[
3{,}1408450 < \pi < 3{,}1428571.
\]  
Aber auch dies ist eben wieder nur eine Näherung. Neben Verfeinerungen der Methode von Archimedes und der Entwicklung anderen Zugänge, wurde die Zahl $\pi$ auf immer mehr Nachkommastellen approximiert. Die Frage nach der Rationalität von $\pi$ drängte sich förmlich auf (Archimedes zweifelte dies bereits an), aber blieb trotz intensiver Untersuchungen lange offen. Es scheint, dass ein erster Beweis, der mittels Kettenbrüchen (Latein: fractio continua) geführt wurde, auf den Mathematiker Johann Heinrich Lambert und die Jahre um 1770 zurückgeht. Er schreibt in seiner Abhandlung \cite{L1770}:
\vskip 2mm
"`\emph{Ich habe aber in vorbemeldeter Abhandlung von Verwandlung der Brüche die andere Fractio continua angegeben, welche nach einem gewissen Gesetze ins Unendliche Fortgeht, und die Hofnung, die Verhältniß des Diameters zum Umkreise durch ganze Zahlen zu bestimmen, ganz benimmt.}"' 
\vskip 2mm

Trotz (oder wegen?) der revolutionären Erkenntnis, dass nie alle
Stellen von $\pi$ bestimmt werden können, hat die Suche nach weiteren Stellen
nie aufgehört, und verwendet heute weit aufwändigere mathematische 
Methoden als der Beweis der Irrationalität.  Zum Zeitpunkt der Entstehung dieses
Artikels waren übrigens die ersten 31{,}4 Billionen Nachkommastellen von $\pi$
bekannt.
Rationale Näherungen von $\pi$, welche es in den "`Mainstream"' geschafft haben, 
sind $3{,}14=\frac{314}{100}$ und $\frac{22}{7}=3{,}142857\dots$(Mathematik-begeisterte begehen deswegen am 14.~März den "`$\pi$-day"' und am 
22.~Juli den "`$\pi$-approximation-day"').

Wir wollen nun mit den Grundlagen fortfahren, um einen modernen Beweis für die Irrationalität zu führen. Zu diesem Zweck beschäftigen wir uns im folgenden Abschnitt mit dem Integral sowie dem Hauptsatz der Differential- und Integralrechnung.

% % % % % % % % % % % % % % % % % % % % % %
\subsection{Das Integral und der Hauptsatz}
% % % % % % % % % % % % % % % % % % % % % %

Das Integral $\int_a^b f(x)\,dx$ über eine nicht-negative, stetige Funktion in den 
Grenzen $a$ und $b$ hat die Interpretation des Inhalts der Fläche
zwischen der $x$-Achse und dem Graphen von $f$ und zwischen $a$ und $b$.
\begin{center}
%% Creator: Inkscape 1.0 (4035a4f, 2020-05-01), www.inkscape.org
%% PDF/EPS/PS + LaTeX output extension by Johan Engelen, 2010
%% Accompanies image file 'def-int.pdf' (pdf, eps, ps)
%%
%% To include the image in your LaTeX document, write
%%   \input{<filename>.pdf_tex}
%%  instead of
%%   \includegraphics{<filename>.pdf}
%% To scale the image, write
%%   \def\svgwidth{<desired width>}
%%   \input{<filename>.pdf_tex}
%%  instead of
%%   \includegraphics[width=<desired width>]{<filename>.pdf}
%%
%% Images with a different path to the parent latex file can
%% be accessed with the `import' package (which may need to be
%% installed) using
%%   \usepackage{import}
%% in the preamble, and then including the image with
%%   \import{<path to file>}{<filename>.pdf_tex}
%% Alternatively, one can specify
%%   \graphicspath{{<path to file>/}}
%% 
%% For more information, please see info/svg-inkscape on CTAN:
%%   http://tug.ctan.org/tex-archive/info/svg-inkscape
%%
\begingroup%
  \makeatletter%
  \providecommand\color[2][]{%
    \errmessage{(Inkscape) Color is used for the text in Inkscape, but the package 'color.sty' is not loaded}%
    \renewcommand\color[2][]{}%
  }%
  \providecommand\transparent[1]{%
    \errmessage{(Inkscape) Transparency is used (non-zero) for the text in Inkscape, but the package 'transparent.sty' is not loaded}%
    \renewcommand\transparent[1]{}%
  }%
  \providecommand\rotatebox[2]{#2}%
  \newcommand*\fsize{\dimexpr\f@size pt\relax}%
  \newcommand*\lineheight[1]{\fontsize{\fsize}{#1\fsize}\selectfont}%
  \ifx\svgwidth\undefined%
    \setlength{\unitlength}{277.5359159bp}%
    \ifx\svgscale\undefined%
      \relax%
    \else%
      \setlength{\unitlength}{\unitlength * \real{\svgscale}}%
    \fi%
  \else%
    \setlength{\unitlength}{\svgwidth}%
  \fi%
  \global\let\svgwidth\undefined%
  \global\let\svgscale\undefined%
  \makeatother%
  \begin{picture}(1,0.51776043)%
    \lineheight{1}%
    \setlength\tabcolsep{0pt}%
    \put(0,0){\includegraphics[width=\unitlength,page=1]{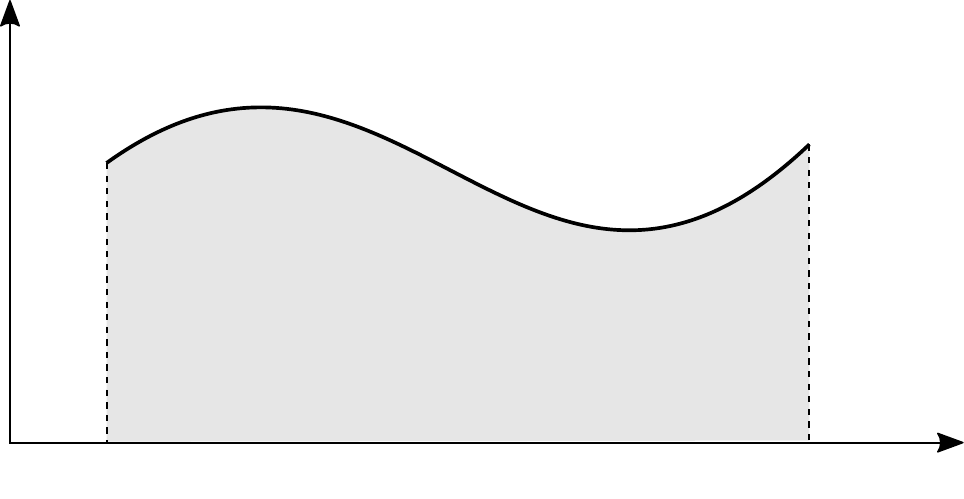}}%
    \put(0.1094558,0.00790196){\makebox(0,0)[t]{\lineheight{1.25}\smash{\begin{tabular}[t]{c}$a$\end{tabular}}}}%
    \put(0.84103897,0.00682707){\makebox(0,0)[t]{\lineheight{1.25}\smash{\begin{tabular}[t]{c}$b$\end{tabular}}}}%
    \put(0.56578627,0.32495474){\makebox(0,0)[t]{\lineheight{1.25}\smash{\begin{tabular}[t]{c}$f$\end{tabular}}}}%
    \put(0.47560645,0.19855188){\makebox(0,0)[t]{\lineheight{1.25}\smash{\begin{tabular}[t]{c}$A=\int_a^b f(x)dx$\end{tabular}}}}%
  \end{picture}%
\endgroup%

\end{center}
Geometrisch ist folgende Eigenschaft des Integrals unmittelbar einleuchtend:
Liegt der Graph von $g$ stets über dem von $f$, so ist
der Inhalt der begrenzten Fläche größer, somit auch das Integral. 
Mathematisch drücken wir dies so aus:
\begin{satz*}[Monotonie des Integrals] 
Sind $f,g\colon [a;b]\to\R$ stetige Funktionen mit der Eigenschaft,
dass für alle $x\in [a;b]$ die Ungleichung $g(x)\le f(x)$ gilt, so ist
\[
\int_a^b g(x)\,dx\le \int_a^b f(x) \,dx\;.
\]
Speziell:  Gilt für Zahlen $c,d\in\R$ die Ungleichung $c\le f(x)\le d$ für alle $x\in\R$,
dann folgt
\[
(b-a)\cdot c\le \int_a^bf(x)\, dx\le (b-a)\cdot d\;.
\]
\end{satz*}
Die folgende Grafik illustriert eine zur Monotonie äquivalente Aussage, nämlich 
die, dass wenn der Graph von $f$ über dem von $g$ liegt, die von
$f$ und $g$ begrenzte Fläche nicht-negativ ist:
\begin{center}
%% Creator: Inkscape 1.0 (4035a4f, 2020-05-01), www.inkscape.org
%% PDF/EPS/PS + LaTeX output extension by Johan Engelen, 2010
%% Accompanies image file 'monotonie-int.pdf' (pdf, eps, ps)
%%
%% To include the image in your LaTeX document, write
%%   \input{<filename>.pdf_tex}
%%  instead of
%%   \includegraphics{<filename>.pdf}
%% To scale the image, write
%%   \def\svgwidth{<desired width>}
%%   \input{<filename>.pdf_tex}
%%  instead of
%%   \includegraphics[width=<desired width>]{<filename>.pdf}
%%
%% Images with a different path to the parent latex file can
%% be accessed with the `import' package (which may need to be
%% installed) using
%%   \usepackage{import}
%% in the preamble, and then including the image with
%%   \import{<path to file>}{<filename>.pdf_tex}
%% Alternatively, one can specify
%%   \graphicspath{{<path to file>/}}
%% 
%% For more information, please see info/svg-inkscape on CTAN:
%%   http://tug.ctan.org/tex-archive/info/svg-inkscape
%%
\begingroup%
  \makeatletter%
  \providecommand\color[2][]{%
    \errmessage{(Inkscape) Color is used for the text in Inkscape, but the package 'color.sty' is not loaded}%
    \renewcommand\color[2][]{}%
  }%
  \providecommand\transparent[1]{%
    \errmessage{(Inkscape) Transparency is used (non-zero) for the text in Inkscape, but the package 'transparent.sty' is not loaded}%
    \renewcommand\transparent[1]{}%
  }%
  \providecommand\rotatebox[2]{#2}%
  \newcommand*\fsize{\dimexpr\f@size pt\relax}%
  \newcommand*\lineheight[1]{\fontsize{\fsize}{#1\fsize}\selectfont}%
  \ifx\svgwidth\undefined%
    \setlength{\unitlength}{281.368651bp}%
    \ifx\svgscale\undefined%
      \relax%
    \else%
      \setlength{\unitlength}{\unitlength * \real{\svgscale}}%
    \fi%
  \else%
    \setlength{\unitlength}{\svgwidth}%
  \fi%
  \global\let\svgwidth\undefined%
  \global\let\svgscale\undefined%
  \makeatother%
  \begin{picture}(1,0.60805956)%
    \lineheight{1}%
    \setlength\tabcolsep{0pt}%
    \put(0,0){\includegraphics[width=\unitlength,page=1]{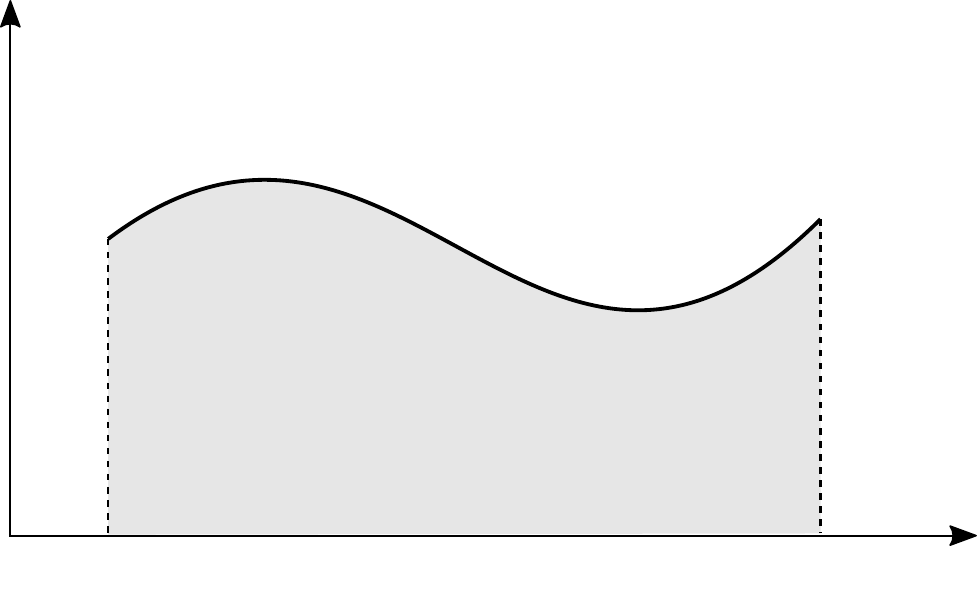}}%
    \put(0.10966278,0.00688429){\makebox(0,0)[t]{\lineheight{1.25}\smash{\begin{tabular}[t]{c}$a$\end{tabular}}}}%
    \put(0.84100625,0.00575918){\makebox(0,0)[t]{\lineheight{1.25}\smash{\begin{tabular}[t]{c}$b$\end{tabular}}}}%
    \put(0.56584374,0.33879612){\makebox(0,0)[t]{\lineheight{1.25}\smash{\begin{tabular}[t]{c}$f$\end{tabular}}}}%
    \put(0,0){\includegraphics[width=\unitlength,page=2]{monotonie-int.pdf}}%
    \put(0.56750843,0.1783407){\makebox(0,0)[t]{\lineheight{1.25}\smash{\begin{tabular}[t]{c}$g$\end{tabular}}}}%
    \put(0.34495652,0.2577913){\makebox(0,0)[t]{\lineheight{1.25}\smash{\begin{tabular}[t]{c}$\int_a^b \big(f(x)-g(x)\big)dx$\end{tabular}}}}%
  \end{picture}%
\endgroup%

\end{center}
Die folgende Grafik veranschaulicht die zweite Aussage im vorangehenden Satz für eine nicht-negative, stetige Funktion, nämlich dass die Fläche unterhalb des Graphen von $f$ in den Grenzen $a$ bis $b$ zwar größer ist als die Fläche des Rechteckes mit den Kantenlängen $(b-a)$ und $c$, aber kleiner ist als die Fläche des Rechteckes mit Kantenlängen $(b-a)$ und $d$.
\begin{center}
%% Creator: Inkscape 1.0 (4035a4f, 2020-05-01), www.inkscape.org
%% PDF/EPS/PS + LaTeX output extension by Johan Engelen, 2010
%% Accompanies image file 'schranken-int.pdf' (pdf, eps, ps)
%%
%% To include the image in your LaTeX document, write
%%   \input{<filename>.pdf_tex}
%%  instead of
%%   \includegraphics{<filename>.pdf}
%% To scale the image, write
%%   \def\svgwidth{<desired width>}
%%   \input{<filename>.pdf_tex}
%%  instead of
%%   \includegraphics[width=<desired width>]{<filename>.pdf}
%%
%% Images with a different path to the parent latex file can
%% be accessed with the `import' package (which may need to be
%% installed) using
%%   \usepackage{import}
%% in the preamble, and then including the image with
%%   \import{<path to file>}{<filename>.pdf_tex}
%% Alternatively, one can specify
%%   \graphicspath{{<path to file>/}}
%% 
%% For more information, please see info/svg-inkscape on CTAN:
%%   http://tug.ctan.org/tex-archive/info/svg-inkscape
%%
\begingroup%
  \makeatletter%
  \providecommand\color[2][]{%
    \errmessage{(Inkscape) Color is used for the text in Inkscape, but the package 'color.sty' is not loaded}%
    \renewcommand\color[2][]{}%
  }%
  \providecommand\transparent[1]{%
    \errmessage{(Inkscape) Transparency is used (non-zero) for the text in Inkscape, but the package 'transparent.sty' is not loaded}%
    \renewcommand\transparent[1]{}%
  }%
  \providecommand\rotatebox[2]{#2}%
  \newcommand*\fsize{\dimexpr\f@size pt\relax}%
  \newcommand*\lineheight[1]{\fontsize{\fsize}{#1\fsize}\selectfont}%
  \ifx\svgwidth\undefined%
    \setlength{\unitlength}{277.5359159bp}%
    \ifx\svgscale\undefined%
      \relax%
    \else%
      \setlength{\unitlength}{\unitlength * \real{\svgscale}}%
    \fi%
  \else%
    \setlength{\unitlength}{\svgwidth}%
  \fi%
  \global\let\svgwidth\undefined%
  \global\let\svgscale\undefined%
  \makeatother%
  \begin{picture}(1,0.51776043)%
    \lineheight{1}%
    \setlength\tabcolsep{0pt}%
    \put(0,0){\includegraphics[width=\unitlength,page=1]{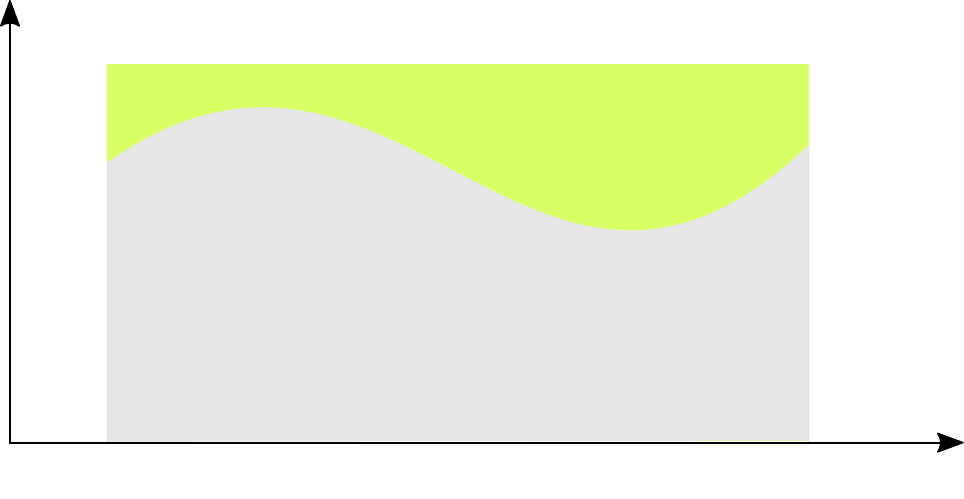}}%
    \put(0.1094558,0.00790196){\makebox(0,0)[t]{\lineheight{1.25}\smash{\begin{tabular}[t]{c}$a$\end{tabular}}}}%
    \put(0.84103897,0.00682707){\makebox(0,0)[t]{\lineheight{1.25}\smash{\begin{tabular}[t]{c}$b$\end{tabular}}}}%
    \put(0.56578627,0.32495474){\makebox(0,0)[t]{\lineheight{1.25}\smash{\begin{tabular}[t]{c}$f$\end{tabular}}}}%
    \put(0.47560645,0.19855188){\makebox(0,0)[t]{\lineheight{1.25}\smash{\begin{tabular}[t]{c}$A=\int_a^b f(x)dx$\end{tabular}}}}%
    \put(0,0){\includegraphics[width=\unitlength,page=2]{schranken-int.pdf}}%
    \put(0.85482483,0.22877378){\makebox(0,0)[lt]{\lineheight{1.25}\smash{\begin{tabular}[t]{l}$c$\end{tabular}}}}%
    \put(0.85332321,0.44256822){\makebox(0,0)[lt]{\lineheight{1.25}\smash{\begin{tabular}[t]{l}$d$\end{tabular}}}}%
  \end{picture}%
\endgroup%

\end{center}

Wir betrachten nun zur stetigen Funktion $f$ auf $[a;b]$ die 
Funktion $F$ auf $[a;b]$ mit 
\[
F(x)=\int_a^x f(t) \,dt,
\]
d.h. die Funktion, die ein Element aus $x\in[a,b]$ auf das Integral über die Funktion $f$ in den Integrationsgrenzen $a$ bis $x$ abbildet. 
\begin{figure}[H]
\begin{center}
%% Creator: Inkscape 1.0 (4035a4f, 2020-05-01), www.inkscape.org
%% PDF/EPS/PS + LaTeX output extension by Johan Engelen, 2010
%% Accompanies image file 'hauptsatz.pdf' (pdf, eps, ps)
%%
%% To include the image in your LaTeX document, write
%%   \input{<filename>.pdf_tex}
%%  instead of
%%   \includegraphics{<filename>.pdf}
%% To scale the image, write
%%   \def\svgwidth{<desired width>}
%%   \input{<filename>.pdf_tex}
%%  instead of
%%   \includegraphics[width=<desired width>]{<filename>.pdf}
%%
%% Images with a different path to the parent latex file can
%% be accessed with the `import' package (which may need to be
%% installed) using
%%   \usepackage{import}
%% in the preamble, and then including the image with
%%   \import{<path to file>}{<filename>.pdf_tex}
%% Alternatively, one can specify
%%   \graphicspath{{<path to file>/}}
%% 
%% For more information, please see info/svg-inkscape on CTAN:
%%   http://tug.ctan.org/tex-archive/info/svg-inkscape
%%
\begingroup%
  \makeatletter%
  \providecommand\color[2][]{%
    \errmessage{(Inkscape) Color is used for the text in Inkscape, but the package 'color.sty' is not loaded}%
    \renewcommand\color[2][]{}%
  }%
  \providecommand\transparent[1]{%
    \errmessage{(Inkscape) Transparency is used (non-zero) for the text in Inkscape, but the package 'transparent.sty' is not loaded}%
    \renewcommand\transparent[1]{}%
  }%
  \providecommand\rotatebox[2]{#2}%
  \newcommand*\fsize{\dimexpr\f@size pt\relax}%
  \newcommand*\lineheight[1]{\fontsize{\fsize}{#1\fsize}\selectfont}%
  \ifx\svgwidth\undefined%
    \setlength{\unitlength}{336.67450441bp}%
    \ifx\svgscale\undefined%
      \relax%
    \else%
      \setlength{\unitlength}{\unitlength * \real{\svgscale}}%
    \fi%
  \else%
    \setlength{\unitlength}{\svgwidth}%
  \fi%
  \global\let\svgwidth\undefined%
  \global\let\svgscale\undefined%
  \makeatother%
  \begin{picture}(1,0.68681967)%
    \lineheight{1}%
    \setlength\tabcolsep{0pt}%
    \put(0,0){\includegraphics[width=\unitlength,page=1]{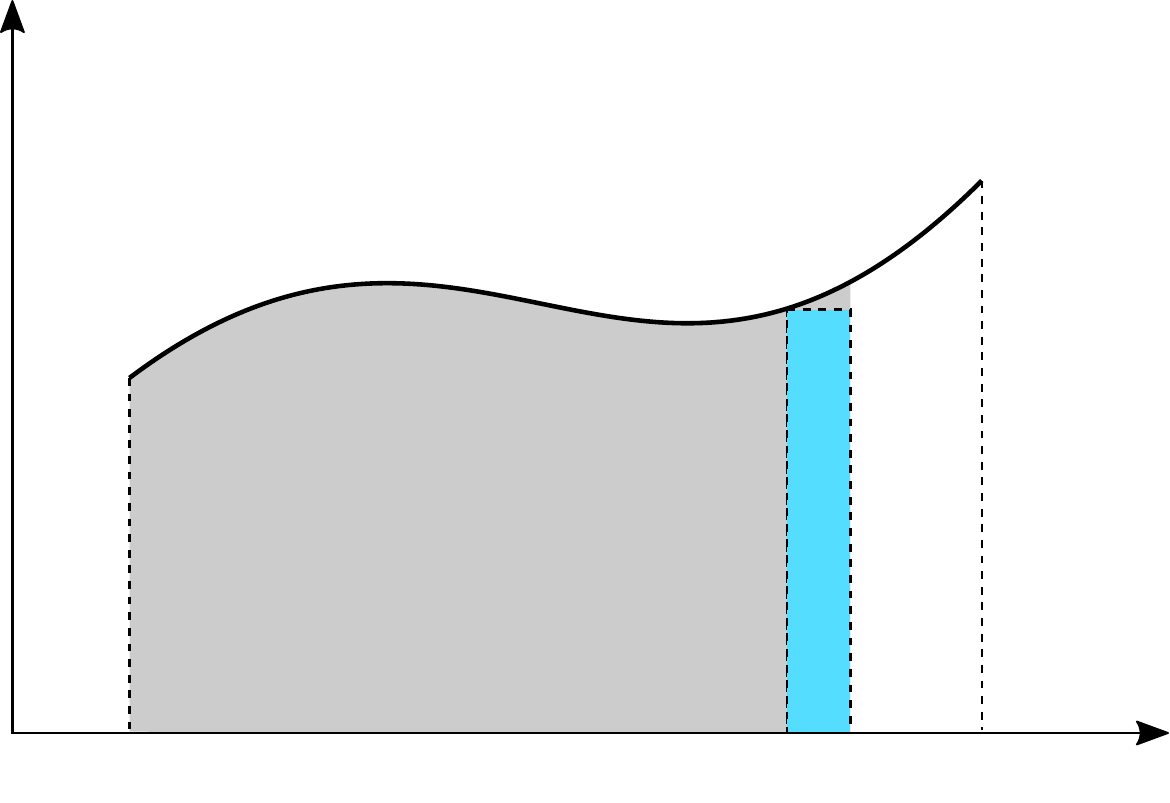}}%
    \put(0.10966277,0.00688433){\makebox(0,0)[t]{\lineheight{1.25}\smash{\begin{tabular}[t]{c}$a$\end{tabular}}}}%
    \put(0.84100623,0.00575917){\makebox(0,0)[t]{\lineheight{1.25}\smash{\begin{tabular}[t]{c}$b$\end{tabular}}}}%
    \put(0.68335569,0.00731889){\makebox(0,0)[rt]{\lineheight{1.25}\smash{\begin{tabular}[t]{r}$x$\end{tabular}}}}%
    \put(0.74528127,0.00619367){\makebox(0,0)[t]{\lineheight{1.25}\smash{\begin{tabular}[t]{c}$x+h$\end{tabular}}}}%
    \put(0.61485221,0.2274125){\makebox(0,0)[t]{\lineheight{1.25}\smash{\begin{tabular}[t]{c}$f(x)$\end{tabular}}}}%
    \put(0,0){\includegraphics[width=\unitlength,page=2]{hauptsatz.pdf}}%
  \end{picture}%
\endgroup%

\caption{Der Hauptsatz  der Differential- und Integralrechnung}\label{fig:hauptsatz}
\end{center}
\end{figure}

Man sieht aus Abbildung \ref{fig:hauptsatz}, dass für (hinreichend kleine) $h\in(0;\infty)$
\[
F(x+h)-F(x)=\int_a^{x+h} f(t)\,dt- \int_a^x f(t)\,dt
\approx h\cdot f(x)\;,
\]
wobei aufgrund der Stetigkeit von $f$ die Näherung\footnote{Das Zeichen $\approx$, welches die Näherung andeutet, verwenden wir an dieser Stelle nur intuitiv, ohne exakt definiert zu haben, was es bedeuten soll. Wir wollen es an dieser Stelle bei folgender Erklärung belassen: der Wert
auf der linken Seite stimmt mit dem auf der rechten Seite bis
auf einen gewissen Fehler überein. Der Fehler verschwindet, wenn $h$ gegen 0 geht.} desto besser ist, je kleiner $h$ ist.
Es gilt also für den Grenzwert des Differenzenquotienten
\[
\lim_{h\to 0}\frac{F(x+h)-F(x)}{h}=f(x)\;,
\]
mit anderen Worten, $F'(x)=f(x)$.\\

Was wir uns soeben anschaulich überlegt haben, ist ein Teil der Aussage des Hauptsatzes der Differential- und Integralrechnung, den wir jetzt formulieren möchten. 

\begin{satz*}[Hauptsatz der Differential- und Integralrechnung]
Es sei $f\colon[a;b]\to\R$ eine stetige Funktion auf dem Intervall $[a;b]$. Dann gilt für alle $x_0\in[a;b]$, dass die Funktion
\[
F:[a;b]\to\R,\qquad F(x)=\int_{x_0}^x f(t) \,dt
\]
differenzierbar ist und für alle $x\in[a;b]$, $F'(x)=f(x)$ gilt. Insbesondere haben wir
\[
\int_{a}^b f(t)\,dt = F(b) - F(a).
\]
\end{satz*}

Bevor wir mit dem Beweis der Irrationalität von $\pi$ fortfahren, wollen wir noch anmerken, dass es natürlich viele Funktionen geben kann, deren Ableitung $f$ ist. Allerdings unterscheiden sich diese nur bis auf eine Konstante. Diese Aussage ist Teil der folgenden Bemerkung.

\begin{rem*}
Ist $G$ eine beliebige Funktion auf $[a;b]$ mit Ableitung $f$, so ist 
also $(F-G)'=F'-G'=f-f=0$, also $F-G$ eine konstante Funktion. 
Es gibt also eine Konstante $C\in\R$ mit $F(x)=G(x)+C$ für alle $x\in[a;b]$.
Es ist damit 
\begin{align*}
\int_a^b f(t)\,dt=F(b)-F(a)=\big(G(b)+C\big)-\big(G(a)+C\big)=G(b)-G(a)\,.
\end{align*}
\end{rem*}

% % % % % % % % % % % % % % % % % %
\section{Die Irrationalität von $\pi$}
% % % % % % % % % % % % % % % % % %

Der Beweis der Irrationalität von $\pi$, den wir hier präsentieren, stammt von Ivan Niven \cite{N1947} und erschien im Jahre 1947 und wird auch im Lehrbuch zur Analysis 1 von Königsberger behandelt \cite{K1990}.

Wenn wir zeigen können, dass $\pi^2\notin \Q$, dann ist auch $\pi\notin \Q$:
Angenommen es gäbe $a,b\in \N$ mit  $\pi=\frac{a}{b}$, dann 
wäre $\pi^2=\frac{a^2}{b^2}$, und somit rational.

\begin{satz*}[Irrationalität von $\pi$]
Es gibt keine natürlichen Zahlen $a,b\in\N$ mit $\pi^2=\frac{a}{b}\,$. 
\end{satz*}

\begin{proof}
Wir nehmen das Gegenteil der Aussage an  und führen diese Annahme auf einen 
Widerspruch. 

Es seien also $a,b\in \N$ mit $\pi^2=\frac{a}{b}$. Wir wählen 
$n\in\N$ groß genug, so dass $\frac{\pi a^n}{n!}<1$. Dass ein solches $n\in\N$ existiert, haben wir bereits im Abschnitt über die Faktorielle bewiesen. Wir definieren nun eine Funktion $f\colon [0,1]\to \R$ durch
\[
f(x):=\frac{1}{n!}x^n\cdot (1-x)^n\;.
\]
Mithilfe des binomischen Lehrsatzes erhalten wir
\begin{align*}
f(x)
&=\frac{1}{n!}x^n\sum_{k=0}^n {n\choose k}(-1)^{k}x^k
=\sum_{k=0}^n \frac{1}{k!(n-k)!}(-1)^{k}x^{n+k}\\
&=\sum_{j=n}^{2n} \frac{1}{(j-n)!(2n-j)!}(-1)^{n-j}x^{j}\;.
\end{align*}
Für alle $\ell\in\N\cup\{0\}$ mit $n\le \ell\le 2 n$, gilt nach unserem Satz über die höhere Ableitung von Potenzfunktionen
\begin{align*}
f^{(\ell)}(x)
&=\sum_{j=\ell}^{2n} \frac{1}{(j-n)!(2n-j)!}(-1)^{n-j}\frac{j!}{(j-\ell)!}x^{j-\ell}\;
\end{align*}
und somit durch Auswertung an der Stelle 0
\begin{align*}
f^{(\ell)}(0)
&= \frac{1}{(\ell-n)!(2n-\ell)!}(-1)^{n-\ell}\frac{\ell!}{0!}
= {n\choose \ell-n}(-1)^{n-\ell}\frac{\ell!}{n!}\;.
\end{align*}
Für alle $\ell\in\N\cup\{0\}$ mit $\ell<n$ oder $\ell>2 n$, gilt automatisch $f^{(\ell)}(0)=0$, 
so dass wir insgesamt für alle $\ell\in \N\cup\{0\}$ die Aussage $f^{(\ell)}(0)\in \Z$ erhalten.
Da für alle $x\in [0;1]$ gilt, dass $f(1-x)=f(x)$, ist auch 
$f^{(\ell)}(1)\in \Z$ für jedes $\ell\in \N\cup\{0\}$.

Wir definieren nun eine weitere Funktion $g\colon [0;1]\to \R$ durch
\begin{align*}
g(x)& := b^n\sum_{k=0}^n (-1)^k\pi^{2n-2k}f^{(2k)}(x) \cr
& = b^n\cdot\big(\pi^{2n}f(x)-\pi^{2n-2}f''(x)
+\pi^{2n-4}f^{(4)}(x)-\ldots+(-1)^n f^{(2n)}(x)\big)\;.
\end{align*}
Da wegen der Annahme $\pi^2 = \frac{a}{b}$ mit $a,b\in\N$ die Zahl $a=b \pi^2$ eine natürliche Zahl ist, sind $g(0)$ und $g(1)$ ganze Zahlen.  
Weiter ergibt sich, wegen der Identität
\begin{align*}
g''(x) + \pi^2g(x) & = b^n\sum_{k=0}^n (-1)^k\pi^{2n-2k} f^{(2k+2)}(x) + \pi^2 b^n \sum_{k=0}^n (-1)^k\pi^{2n-2k} f^{(2k)}(x)\cr
& = \pi^2b^n(-1)^{0}\pi^{2n-2\cdot0}f^{(0)}(x),
\end{align*}
dass mit Produkt- und  Kettenregel für das Differenzieren die Gleichheit
\begin{align*}
\Big(g'(x)\sin(\pi x)-\pi g(x) \cos(\pi x)\Big)'
&=\big(g''(x)+\pi^2 g(x)\big)\sin(\pi x)\\
&=b^n\pi^{2n+2}f(x)\sin(\pi x)\\
&=a^n\pi^{2}f(x)\sin(\pi x)\;
\end{align*}
gilt. \footnote{Und genau hier kommt die Magie ins Spiel! Die Wahl der Funktion $g$ ist so geschickt, dass sie wie auf zauberhafte Weise genau so mit den trigonometrischen Funktionen und der Zahl $\pi$ interagiert, wie es benötigt wird.}

Mit dem Hauptsatz der Differential- und Integralrechnung sowie der sich anschließenden Bemerkung, ist also 
\begin{align*}
I&:=\pi a^n\int_0^1f(x)\sin(\pi x)\,dx\\
&=\frac{1}{\pi}\big(g'(x)\sin(\pi x)-\pi g(x) \cos(\pi x)\big)\Big|_0^1
=g(0)+g(1)\;.
\end{align*}
Damit ist nach den vorangehenden Überlegungen $I$ eine ganze Zahl. Andererseits gilt für alle $x\in [0;1]$,
dass $f(x)<\frac{1}{n!}$ und $0\le \sin(\pi x)\le 1$. Aus der
Monotonie des Integrals ergibt sich damit die Abschätzung 
\[
0<I<\frac{\pi a^n}{n!}<1\;.
\]
Dies steht im Widerspruch zu unserer früheren Feststellung, wonach $I\in \Z$. Damit kann $\pi^2$ nicht rational sein.
\end{proof}

Neben elementaren Rechenregeln und den in den einführenden Kapiteln beschriebenen Methoden wurden im Beweis für die Irrationalität von $\pi$
lediglich die Produktregel für das Differenzieren sowie die Kettenregel in ihrer einfachsten Form verwendet.

% % % % % % % % % % % % % % % % % % % %
\section{Die Irrationalität von $e$}
% % % % % % % % % % % % % % % % % % % %

Wir wollen nun zum Abschluss noch zeigen, dass auch 
die Eulersche Zahl $e$ irrational ist. Dazu benötigen wir eine Darstellung von $e$, oder genauer gesagt, von der Exponentialfunktion $x\mapsto e^x$, die typischerweise nicht in der Schule behandelt wird. Während diese Funktion in der Schule über Folgen eingeführt wird, nämlich durch
\[
e^x := \lim_{n\to\infty}\Big(1+\frac{x}{n}\Big)^n,\qquad x\in\R,
\]
benötigen wir die Entwicklung dieser Funktion als Potenzreihe. In dieser Darstellung haben wir, wie wir gleich beweisen werden,
\[
e^x = \sum_{k=0}^\infty \frac{x^k}{k!}.
\] 
Dieser Formel möchten wir noch ein Zitat des bekannten Mathematikers Walter Rudin \cite{R1966} anschließen:

\begin{center}
\includegraphics[width=13cm]{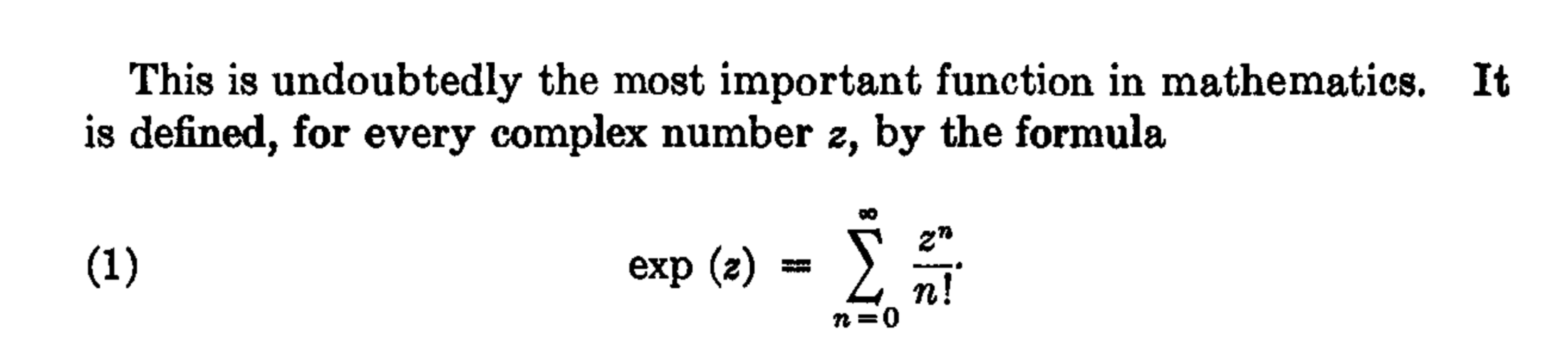}
\end{center}

Entsprechend der Jagd nach immer mehr Nachkommastellen der Zahl $\pi$, hat auch die Suche nach weiteren Stellen bei der Eulerschen Zahl $e$
nie aufgehört. Auch hier kommen heute aufwändige mathematische Methoden zum Einsatz. Zum Zeitpunkt der Entstehung dieses
Artikels waren übrigens die ersten $8$ Billionen Nachkommastellen von $e$
bekannt.

Wir zeigen nun, wie sich aus dem binomischen Lehrsatz die Potenzreihendarstellung für $e$ ableiten lässt.

\begin{satz*} Für alle $x\in(0;\infty)$ ist 
\[
e^x:=\lim_{n\to\infty}\big(1+\frac{x}{n}\big)^n
=\lim_{n\to\infty} 1+x+\frac{x^2}{2!}+\ldots+\frac{x^n}{n!}.
\]
\end{satz*}

\begin{proof}
Für alle $n\in \N$ und alle $k\in\N\cup\{0\}$ mit $0\le k \le n $ ist
\begin{align*}
{n\choose k}\frac{1}{n^k}=\frac{1}{k!}\frac{1\cdot 2\cdot \ldots\cdot (n-k)}{1\cdot 2\cdot \ldots\cdot (n-k)}\frac{(n-k+1)\cdot \ldots\cdot n}{n\cdot\ldots\cdot n}\;.
\end{align*}
Also ergibt sich aus dem binomischen Lehrsatz, dass
\begin{align*}
\big(1+\frac{x}{n}\big)^n
&=1+{n\choose 1}\frac{x}{n}+{n\choose 2}\frac{x^2}{n^2}
+\ldots+{n\choose n-1}\frac{x^{n-1}}{n^{n-1}}+\frac{x^n}{n^n}\\
&<1+x+\frac{x^2}{2!}+\ldots+\frac{x^n}{n!}\;.
\end{align*}
Umgekehrt ist für alle $n,m\in\N$ und alle $k\in\N\cup\{0\}$ mit $0\le k\le n+m$
\[
{m+n\choose k}\frac{1}{(m+n)^k}
=\frac{1}{k!}\frac{1\cdot 2\cdot \ldots\cdot m}{1\cdot 2\cdot \ldots\cdot m}\frac{(m+1)\cdot \ldots\cdot (m+n)}{(m+n)\cdot\ldots\cdot (m+n)}
\stackrel{m\to\infty}{\longrightarrow}\frac{1}{k!}\;.
\]
Somit erhalten wir
\begin{align*}
\big(1+\frac{x}{n+m}\big)^{n+m}
&>1+{n+m\choose 1}\frac{x}{n+m}+\ldots+{n+m\choose n}\frac{x^{n}}{(m+n)^n}\\
&\stackrel{m\to\infty}{\longrightarrow}1+x+\frac{x^2}{2!}+\ldots+\frac{x^n}{n!}\;.
\end{align*}
Dies beweist die gewünschte Darstellung.
\end{proof}

Kommen wir nun zum Beweis der Irrationalität der Eulerschen Zahl, der lediglich auf der Potenzreihendarstellung beruht und einige elementare Abschätzungen verwendet.

\begin{satz*}[Irrationalität von $e$]
Es gibt keine natürlichen Zahlen $a,b\in\N$ mit $e=\frac{a}{b}\,$. 
\end{satz*}

\begin{proof}
Wir nehmen das Gegenteil der Aussage an  und führen diese Annahme auf einen 
Widerspruch. 

Es seien also $a,b\in \N$ mit $e=\frac{a}{b}$. Wir wählen 
$n\in\N$ groß genug, so dass $n!\cdot\frac{a}{b}\in\N$ (man könnte zum Beispiel $n=b$ wählen). Dann gilt

\begin{align*}
e
&=\sum_{k=0}^\infty \frac{1}{k!}
=\sum_{k=0}^{n-1} \frac{1}{k!}
+\sum_{k=n}^\infty \frac{1}{k!}
<\sum_{k=0}^{n-1} \frac{1}{k!}
+\sum_{k=n}^\infty \frac{1}{n!\cdot (n+1)^{k-n}}\\
&=\sum_{k=0}^{n-1} \frac{1}{k!}
+\frac{1}{n!}\sum_{k=n}^\infty \Big(\frac{1}{n+1}\Big)^{k-n}
\stackrel{(*)}{=}\sum_{k=0}^{n-1} \frac{1}{k!} +\frac{1}{n!}\cdot\frac{1}{1-\frac{1}{n+1}}\\
&=\sum_{k=0}^{n-1} \frac{1}{k!} +\frac{1}{n!}\cdot\frac{n+1}{n}
=\sum_{k=0}^{n} \frac{1}{k!} +\frac{1}{n!}\cdot\frac{1}{n},
\end{align*}
wo wir für $(*)$ die geometrische Summenformel verwendet haben. Es ist also 
\begin{align*}
0<e-\sum_{k=0}^{n} \frac{1}{k!}<\frac{1}{n!}\cdot\frac{1}{n}\;.
\end{align*}
Wir haben $n!\cdot\frac{a}{b}\in\N$ und $\sum_{k=0}^{n}\frac{n!}{k!}\in\N$, und
damit $M:=n!\cdot\frac{a}{b}-\sum_{k=0}^{n}\frac{n!}{k!}\in \Z$.
Andererseits ist \begin{align*}
0< M=n!\cdot\Big(e-\sum_{k=0}^{n} \frac{1}{k!}\Big)<\frac{1}{n}\le 1\;,
\end{align*}
ein Widerspruch zu $M\in \Z$.
\end{proof}

% % % % % % % % % % % % % % % % % % % % % % % % % % % %
\section{Abschlussbemerkungen}
% % % % % % % % % % % % % % % % % % % % % % % % % % % %

Die Frage nach der Irrationalität von $e$ und $\pi$ ist nicht nur eine ganz natürliche, sondern auch eine selbst für jüngere Schülerinnen und Schüler verständliche Frage, der man sich etwa im Fall von $\pi$ durch tatsächliches Abmessen von Kreisen mittels Maßbändern oder durch einfache Approximation mit Vielecken nähern kann. Durch das breite Spektrum an Anspruchsniveaus eignet sich das Thema in besonderem Maße sowohl für den Einsatz in der Schule als auch an der Universität. Dass zu einem vollständigen Beweis lediglich elementare Methoden der Analysis benötigt werden (die sich in vielen Fällen auch geometrische motivieren lassen), haben wir hoffentlich überzeugend dargelegt. Des Weiteren lässt der Grad der Genauigkeit im Beweis genügend Spielraum das Anspruchsniveau entsprechend anzupassen. 

Ob man die Thematik als eine elegante und unerwartete Anwendung für den Hauptsatz der Differential- und Integralrechnung behandelt, im Rahmen einer Mathematik AG oder eines Projektes, oder gar in Kombination mit dem Informatikunterricht (etwa numerische Näherungen aus der Potenzreihenentwicklung von $e$ oder den Kettenbrüchen), wir wünschen viel Freude und Erfolg bei der Umsetzung.

% % % % % % % % % % % % % %
\subsection*{Danksagung}
% % % % % % % % % % % % % %
Gunther Leobacher wird über das Projekt F5508-N26 des österreichischen Wissenschaftsfond (FWF) gefördert, welches ein Teilprojekt des Spezialforschungsbereichs "`Quasi-Monte Carlo Methoden: Theorie und Anwendungen"' ist.
Joscha Prochno wird durch den österreichischen Wissenschaftsfond (FWF) über das Projekt P32405 "`Asymptotische Geometrische Analysis und Anwendungen"' gefördert. 

Wir danken unserem Kollegen Michael Schmitz (Flensburg) für hilfreiche Anregungen und Kommentare zu dieser Arbeit.
 
%\bibliographystyle{plain}
%\bibliography{irrat}

\end{document}